\def\arxiv{1}
\if\arxiv1
\title{ \bf%
Some results on exponential synchronization of nonlinear systems (long version)
\thanks{This work was supported by the ANR Project LIMICOS contract number 12BS0300501}}
\else
\title{ \bf%
Some results on exponential synchronization of nonlinear systems
\thanks{This work was supported by the ANR Project LIMICOS contract number 12BS0300501}}
\fi

\author{Vincent Andrieu\thanks{V. Andrieu is with 
Universit\'e Lyon 1, Villeurbanne, France -- CNRS, UMR 5007, LAGEP, France. { vincent.andrieu@gmail.com}
}, Bayu Jayawardhana\thanks{B.~Jayawardhana is with ENTEG, Faculty of Mathematics and Natural Sciences, University of Groningen, the
Netherlands,
        { bayujw@ieee.org, b.jayawardhana@rug.nl}}, Sophie Tarbouriech\thanks{S.~Tarbouriech is with LAAS-CNRS, Universit\'e de Toulouse, CNRS, Toulouse, France. tarbour@laas.fr.}
}

\documentclass[journal]{IEEEtran}
\IEEEoverridecommandlockouts

\usepackage{amsfonts,latexsym,amsmath,amssymb,mathtools}
\usepackage[mathscr]{eucal}
\usepackage{graphicx,color}
\usepackage{comment}
\usepackage{hyperref}
\usepackage{dsfont}

\newtheorem{theorem}{Theorem}
\newtheorem{lemma}{Lemma}
\newtheorem{proposition}{Proposition}

\newtheorem{definition}{Definition}
\newtheorem{assumption}{Assumption}
\newtheorem{remark}{Remark}

\catcode`\@=11
\def\downparenfill{$\m@th\braceld\leaders\vrule\hfill\bracerd$}
\def\overparen#1{\mathop{\vbox{\ialign{##\crcr\crcr
\noalign{\kern0.4ex}
\downparenfill\crcr\noalign{\kern0.4ex\nointerlineskip}
$\hfil\displaystyle{#1}\hfil$\crcr}}}\limits}
\catcode`\@=12

\def\RR{{\mathbb R}}    

\def\CR{\mathcal C}
\DeclareMathOperator{\Id}{I}

\def\de{{\widetilde e}}     
\def\dx{z} 
\def\dX {Z}
\def\dE{{\widetilde E}}     

\def\du{{\tilde u}}     

\def\dk{{\widetilde k}}     

\def\DR{\mathcal D} 		

\def\dz{\tilde z}

\def\bxa{z_a}
\def\bxb{z_b}
\def\bx{z}
\def\dna{n_a}
\def\rhoa{\rho_a}
\def\rhob{\rho_b}
\def\fa{f_a}
\def\fb{f_b}
\def\ga{g_a}
\def\gb{g_b}
\def\PRa{P_a}
\def\PRb{P_b}
\def\Sa{S_a}

\def\qa{q_a}
\def\qb{q_b}
\def\alphaa{\alpha_a}
\def\alphab{\alpha_b}
\def\Ca{\mathcal C_a}
\def\Cb{\mathcal C_b}
\def\va{v_a}
\def\vb{v_b}
\def\Ua{U_a}
\def\Ub{U_b}
\def\Qa{Q_a}
\def\Qb{Q_b}
\def\dG{\bar G}

\def\PR{{P}}

\def \kell{{\ell}}      

\def \La { L_{11}}
\def \Lb {L_{1,2:N}}
\def \Lc {L_{2:N,2:N}}

\newcommand{\bbm}[1]{\left[\begin{matrix} #1 \end{matrix}\right]}

\def\der{{\mathfrak{d}}}

\def\startmodifVA{\color{black}}
\def\stopmodif{\color{black}}

\def\stopmodif{\color{black}}

\begin{document}
\maketitle

\begin{abstract}
Based on recent works on transverse exponential stability,  we establish  some necessary and sufficient conditions for the existence of a (locally) exponential synchronizing control law. We show that the existence of a structured synchronizer is equivalent to the existence of a stabilizer for the individual linearized systems (on the synchronization manifold) by a linear state feedback. This, in turn, is also equivalent to the existence of a symmetric covariant tensor field, which satisfies  a Control Matrix Function   inequality. Based on this  result,  we provide the construction of such synchronizer  via  backstepping approaches. 
In some particular cases, we show how global exponential synchronization may be obtained.
\end{abstract}

\section{Introduction}
\label{Sec_Intro}

Controlled synchronization, as a coordinated control problem of a group of autonomous systems, has been regarded as one of important group behaviors. It has found its relevance in many engineering applications, such as, the distributed control of (mobile) robotic systems, the control and reconfiguration of devices in the context of internet-of-things, and the synchronization of autonomous vehicles (see, for example, \cite{olfatiIEEEproc2007}).   

For linear systems, the solvability of this problem and, as well as, the design of controller, have been thoroughly studied in literature. To name a few, we refer to the classical work on the nonlinear Goodwin oscillators \cite{Goodwin_1965}, to the synchronization of linear systems in \cite{WielandEtAl_Aut_11IntModSyn,ScardoviSepulchre_Aut_09SynchNetw} and to the recent works in nonlinear systems \cite{SarletteSepulchreLeonard_Aut_09,DePersisJayawardhana_TCNS_14,
DePersisJayawardhana_CDC_12,DePersisJayawardhana_SICON_12, ScardoviArcakSontag_10}. For linear systems, the solvability of synchronization problem reduces to the solvability of stabilization of individual systems by either an output or state feedback. It has recently been established in \cite{WielandEtAl_Aut_11IntModSyn} that for linear systems, the solvability of the output synchronization problem is equivalent to the existence of an internal model, which is a well-known concept in the output regulation theory. 

 The generalization of these results to the nonlinear setting has appeared in the literature (see, for example, \cite{Pavlov_Book_2005,Chopra_TAC_2012,Pogromsky_TCS_2001,Isidori_TAC_2014,SarletteSepulchreLeonard_Aut_09,DePersisJayawardhana_TCNS_14,
DePersisJayawardhana_CDC_12,DePersisJayawardhana_SICON_12, ScardoviArcakSontag_10,HamadehStanSepulchreGoncalves_TAC_12}). In these works, the synchronization of nonlinear systems with a fixed network topology can be solved under various different sufficient conditions. 

For instance, the application of passivity theory plays a key role in \cite{Chopra_TAC_2012,Pogromsky_TCS_2001,SarletteSepulchreLeonard_Aut_09,DePersisJayawardhana_SICON_12, ScardoviArcakSontag_10,HamadehStanSepulchreGoncalves_TAC_12}. By using the input/output passivity property, the synchronization control law in these works can simply be given by the relative output measurement. Another approach for synchronizing nonlinear systems is by using output regulation theory as pursued in \cite{Isidori_TAC_2014, DePersisJayawardhana_TCNS_14,Pavlov_Book_2005}. In these papers, the synchronization problem is reformulated as an output regulation problem where the output of each system has to track an exogeneous signal driven by a common exosystem and the resulting synchronization control law is again given by relative output measurement. Lastly, another synchronization approach that has gained interest in recent years is via incremental stability \cite{Angeli_TAC_09_FurthReslyapIncStab} or other related notions, such as, convergent systems \cite{Pavlov_Book_2005}. If we restrict ourselves to the class of incremental ISS, as discussed in \cite{Angeli_TAC_09_FurthReslyapIncStab}, the synchronizer can again be based on the relative output/state measurement. 


Despite assuming a fixed network topology,  necessary and sufficient condition for the solvability of synchronization problem of nonlinear systems is not yet established. Therefore, one of our main contributions of this paper is the characterization of controlled synchronization for general nonlinear systems with fixed network topology.  Using recent results on the transverse exponential contraction, we establish some necessary and sufficient conditions for the solvability of a (locally) exponential synchronization. It extends the work in \cite{AndrieuJayawardhanaPraly_CDC_13} where only two interconnected systems are discussed. 
We show that a necessary condition for achieving synchronization is the existence of  
a symmetric covariant tensor field of order two  whose Lie derivative has to satisfy a Control Matrix Function (\textsf{CMF}) inequality, which is similar to the Control Lyapunov Function and detailed later in Section III.

\if \arxiv1
This paper  extends our preliminary work presented in   \cite{AndrieuJayawardhanaTarbouriech_CDC_Synchro}. In particular, we provide detailed proofs for all main results (which were exempted from the aforementioned paper) and additionally, we present the backstepping approach that allows us to construct a \textsf{CMF}-based synchronizer, as well as, the extension of the local synchronization result to the global one for a specific case.  
\else
This paper  extends our preliminary work presented in   \cite{AndrieuJayawardhanaTarbouriech_CDC_Synchro}.
In particular we improve some results by relaxing some conditions (see the necessary condition section). 
Additionally, we present the backstepping approach that allows us to construct a \textsf{CMF}-based synchronizer, as well as, the extension of the local synchronization result to the global one for a specific case. 
Note that  all proofs are given in the long version of this paper in \cite{AndrieuJayawardhanaTarbouriech_arXiv_16_Synchro}.
\fi

The paper is organized as follows.  We present the problem formulation of synchronization in Section \ref{Sec_ProbDef}. In Section \ref{Sec_NecSufCond}, we present our first main results on  necessary conditions to the solvability of the  synchronization problem.  
Some sufficient conditions for local or global synchronization are given in Section \ref{Sec_subSuffCond}.
A constructive synchronizer design is presented in Section \ref{Sec_BackStep}, where a backstepping procedure is given for designing a \textsf{CMF}-based synchronizing control law.  
\par\vspace{0.2em}
\noindent{\bf Notation.} The vector of all ones with a dimension $N$ is denoted by $\mathds{1}_N$. We denote the identity matrix of dimension $n$ by $\Id_n$ or $\Id$ when no confusion is possible. 
Given $M_1,\dots, M_N$ square matrices, $\texttt{diag}\{M_1,\dots, M_N\}$ is the matrix defined as
$$\texttt{diag}\{M_1,\dots, M_N\}=\begin{bmatrix}
M_1&&\\
&\ddots&\\
&&M_N
\end{bmatrix}.$$
Given a vector field $f$ on $\RR^n$ and a covariant two tensor $P:\RR^n\rightarrow\RR^{m\times m}$, 
$P$ is said to have a derivative along $f$ denoted $\der_{f}P$
if the following limit exists
\begin{equation}
\label{LP11}
\der_{f} P(\dx)\;=\;
 \lim_{h\to 0}
\frac{ P( \dX( \dx,h))- P(\dx)}{h}\ ,
\end{equation}
where $\dX(\dx,\cdot)$ is the flow of the vector field $f$ with an initial state $\dx$ in $\RR^n$.
In that case and, when $m=n$ and $f$ is $C^1$
 $L_fP$ is the Lie derivative of the tensor along $f$ which is defined as
\startmodifVA
\begin{equation}\label{def_LiederTens}
L_fP(\dx) = \der_{f} P(\dx) + P(\dx)\frac{\partial f}{\partial \dx}(\dx) + \frac{\partial f}{\partial \dx}(\dx)^\top P(\dx)\ .
\end{equation}
\stopmodif

\section{Problem definition}
\label{Sec_ProbDef}
\subsection{System description and communication topology}
\label{Sec_Prelim}


In this note, we consider the problem of synchronizing $N$ identical nonlinear systems with $N\geq 2$. For every $i=1,\ldots,N$, the $i$-th system $\Sigma_i$ is described by  
\begin{equation}\label{eq_Syst}
\dot x_i=f(x_i)+g(x_i)u_i\ , \ i=1,\dots, N\,
\end{equation}
where $x_i\in \RR^n$, $u_i\in \RR^p$ and the functions $f$ and $g$ are assumed to be $C^2$. In this setting, all systems has the same drift vector field $f$ and the same control vector field $g:\RR^n\rightarrow\RR^{n\times p}$, but not the same controls in $\RR^p$. For simplicity of notation, we denote the complete state variables by $x = \bbm{x_1^\top\ldots x_N^\top}^\top$ in $\RR^{Nn}$. 
 
The synchronization manifold $\DR$, where the state variables of different systems agree with each other, is defined by 
$$
\DR = \{(x_1,\ldots, x_N)\in\RR^{Nn}\ | \ x_1=x_2= \dots =x_N\}.
$$
For every $x$ in $\RR^{Nn}$, we denote the Euclidean distance to the set $\DR$ by $|x|_\DR$.

The communication graph $\mathcal G$, which is used for synchronizing the state through distributed control $u_i$, $i=1, \ldots, N$, is assumed to be an undirected graph and is defined by $\mathcal G = (\mathcal V, \mathcal E)$, where $\mathcal V$ is the set of $N$ nodes (where the $i$-th node is associated to the system $\Sigma_i$) and $\mathcal E\subset \mathcal V\times \mathcal V$ is a set of $M$ edges that define the pairs of communicating systems. Moreover we assume that the graph $\mathcal G$ is connected. 

Let us, for every edge $k$ in $\mathcal G$ connecting node $i$ to node $j$, label one end (e.g., the node $i$) by a positive sign and the other end (e.g., the node $j$) by a negative sign. 
The incidence matrix $D$ that corresponds to $\mathcal G$ is an $N\times M$ matrix such that 
\[
d_{i,k} = \left\{\begin{array}{ll}+1 & \text{if node }i\text{ is the positive end of edge }k \\ -1 & \text{if node }i\text{ is the negative end of edge }k\\0 & \text{otherwise}\end{array}\right.
\]
Using $D$, the Laplacian matrix $L$ can be given by $L=DD^\top$ whose kernel, by the connectedness of $\mathcal G$, is spanned by $\mathds{1}_{N}$. 
\if \arxiv1
We will need the following lemma on the property of $L$ in some results. 
\vspace{0.1cm}

\begin{lemma}\label{lem_Laplacian}
Let $L= \begin{bmatrix}
\La & \Lb\\
\Lb^\top & \Lc
\end{bmatrix}$ be a non-zero balanced Laplacian matrix associated to an undirected graph $\mathcal G$ where $\La$ is a scalar. Then, the eigenvalues of the $(N-1)\times(N-1)$ matrix $\bar L :=\Lc-\mathds{1}_{N-1}\Lb$ are the same as the non-zero eigenvalues of $L$ with the same multiplicity. Moreover, if the graph is connected then $-\bar L$ is Hurwitz. 
\end{lemma}
The proof of Lemma \ref{lem_Laplacian} can be found in Appendix \ref{sec_ProofLemmaLap}. 
\fi
\subsection{Synchronization problem formulation}
Using the description of the interconnected systems via $\mathcal G$, the state synchronization control problem is defined as follows. 

\begin{definition}\label{Def1}
The control laws 
$u_i=\phi_i(x)$, $i=1\dots, N$ solve the \emph{local uniform exponential synchronization} problem for (\ref{eq_Syst}) if the following conditions hold: 
\begin{enumerate}
\item
For all non-communicating pair $(i,j)$ (i.e., $(i,j)\notin \mathcal E$),
$$\frac{\partial \phi_i}{\partial x_j}(x) = \frac{\partial \phi_j}{\partial x_i}(x) = 0\ , \ \forall x\in \RR^{Nn}; $$
\item For all $x\in \DR$, $\phi(x)=0$ (i.e., $\phi$ is zero on $\DR$); and
\item The manifold $\DR$ of the closed-loop system 
\begin{equation}\label{eq_Syst_Sync}
\dot x_i =f(x_i)+g(x_i)\phi_i(x),\ \ i=1,\dots,N
\end{equation}
is uniformly exponentially stable, i.e., there exist positive constants $r$, $k$ and $\lambda>0$ such that for all $x$ in $\RR^{Nn}$ satisfying $|x|_\DR<r$,  
\begin{equation}\label{eq_Beta_exp}
|X(x,t)|_\DR  
\leq  k\exp(-\lambda t)\,  |x|_\DR,
\end{equation}
where $X(x,t)$ denotes the solution initiated from $x$, holds for all $t$ in the time domain of existence of solution. 
\end{enumerate}
When $r=\infty$, it is called the {\em global uniform exponential synchronization}  problem.   \hfill $\triangle$
\end{definition}
\vspace{0.1cm}

In this definition, the condition 1) implies that the solution $u_i$ is a distributed control law that requires only a local state measurement from its neighbors in the graph $\mathcal G$. 

An important feature of our study is that we focus on exponential stabilization of the synchronizing manifold. This allows us to rely on the study developed in \cite{AndrieuJayawardhanaPraly_CDC_13} (or \cite{AndrieuJayawardhanaPraly_TAC_TransExpStab}) in which an infinitesimal characterization of exponential stability of a transverse manifold is given. As it will be shown in the following section this allows us to formalize some necessary and sufficient conditions in terms of matrix functions ensuring the existence of a synchronizing control law.

\section{Necessary conditions}
\label{Sec_NecSufCond}
\subsection{Infinitesimal stabilizability conditions}\label{subsection_necessary}
In \cite{AndrieuJayawardhanaPraly_CDC_13}, a first attempt has been made to give necessary conditions for the existence of an exponentially synchronizing control law for only two agents.
In \cite{AndrieuJayawardhanaPraly_TAC_TransExpStab}, the same problem has been addressed for $N$ agents but without any communication constraints (all agents can communicate with all others).
In both cases, it is shown that assuming some bounds on derivatives of the vector fields and assuming that the synchronizing control law is invariant by permutation of agents,
 the following two properties are necessary conditions. 
\begin{enumerate}
\item[\textsf{IS}] \textbf{Infinitesimal stabilizability.} The couple $(f,g)$ is such that the $n$-dimensional manifold $\{\dz=0\}$ of the transversally linear system
\begin{subequations}
\begin{align}
\label{eq_SysTran1} \dot \dz &= \frac{\partial f}{\partial \dx}(\dx)\dz + g(\dx)\du\\
\label{eq_SysTran2} \dot \dx &= f(\dx)
\end{align}
\end{subequations}
with $\dz$ in $\RR^n$ and $\dx$ in $\RR^n$ is stabilizable by a state feedback that is linear in $\dz$ (i.e., $\du=h(\dx)\dz$ for some function $h:\RR^n\to \RR^{p\times n}$). 
\item[\textsf{CMF}]
\textbf{Control Matrix Function.}
For all positive definite matrix $Q\in\RR^{n\times n}$, there exist a continuous
function $\PR :\RR^{n}\rightarrow\RR^{n\times n}$, which values are symmetric positive definite matrices
and strictly positive real numbers $\underline{p}$ and $\overline{p}$
such that 
\begin{equation}
\label{LP10}
\underline{p}\Id_n\leq P(\dx)\leq \overline{p}\Id_n
\end{equation} 
holds for all $\dx\in\RR^n$,  
and 
the inequality  (see (\ref{LP11}) and (\ref{def_LiederTens})) 
\begin{equation}\label{eq_TensorDerivSynch}
v^\top L_f\PR (\dx)v \leq  -v^\top Qv\ 
\end{equation}
holds for all $(v,\dx)$ in $\RR^n\times\RR^n$ satisfying $v^\top\PR(\dx)g(\dx)=0$. 
\end{enumerate}

An important feature of properties \textsf{IS} and \textsf{CMF} comes from the fact that they are properties of each individual agent, independent of the network topology. 
The first one is a local stabilizability property. 
The second one establishes that there exists a symmetric covariant tensor field of order two denoted by $\PR$ whose Lie derivative satisfies a certain inequality in some specific directions. 
This type of condition can be related to the notion of control Lyapunov function, which is a characterization of stabilizability as studied by Artstein in \cite{Artstein_NLA_83} or Sontag in \cite{Sontag_SICON_83}.
This property can be regarded as an Artstein like condition.
The dual of the \textsf{CMF} property has been thoroughly studied in  \cite{SanfelicePraly_TAC_12} when dealing with an observer design (\cite[Eq. (8)]{SanfelicePraly_TAC_12}, see also \cite{AndrieuJayawardhanaPraly_CDC_13} or \cite{AndrieuBesanconSerres_CDC_13_NecObs}).

\subsection{Necessity of \textsf{IS} and \textsf{CMF} for exponential synchronization}

We show that properties \textsf{IS} and \textsf{CMF} are still necessary conditions if one considers a network of agents with a communication graph $\mathcal G$ as given in \ref{Sec_Prelim}. 
Hence, as this is already the case for linear system, we recover the paradigm, which establishes that a necessary condition for synchronization is a stabilizability property for each individual agent.  
\begin{theorem}
\label{theo_Necessary}
Consider the interconnected systems in (\ref{eq_Syst}) with the communication graph $\mathcal G$ and assume  that there exists a control law $u=\phi(x)$ where $\phi(x)=\bbm{\phi_1^\top(x) & \dots & \phi_N^\top(x)}^\top$ in $\RR^{Np}$ that solves the local uniform exponential synchronization for (\ref{eq_Syst}).  
 Assume moreover that $g$ is bounded, $f$, $g$ and the $\phi_i$'s have bounded first and second derivatives and the closed-loop system is complete. 
 Then properties \textsf{IS} and \textsf{CMF} hold.
\end{theorem}
Note that this theorem  is a refinement of the result which is written in \cite{AndrieuJayawardhanaTarbouriech_CDC_Synchro} since we have removed an assumption related to the structure of the control law.

\if\arxiv0
The proof of this result can be found in  \cite{AndrieuJayawardhanaTarbouriech_arXiv_16_Synchro}. 
\else
\subsection{Proof of Theorem \ref{theo_Necessary}}

\begin{proof}
The first part of the proof is to show that the synchronizing manifold satisfies a transverse uniform exponential stability property. This allows us to use tools developed in \cite{AndrieuJayawardhanaPraly_TAC_TransExpStab} and show a stabilizability property for an $Nn$-dimensional. Employing some kind of Lyapunov projection, we are able to obtain the stabilizability properties  for the $n$-dimensional transversally linear system (\ref{eq_SysTran1}).

Let $e= \bbm{e_2^\top & e_3^\top & \dots & e_N^\top}^\top$ with $e_i=x_i-x_1$, $i=2, \dots N$, and $z=x_1$.
The closed-loop system (\ref{eq_Syst}) with the control law $\phi$ is given by
\begin{equation}\label{eq_System}
\dot e = F(e,z)\ ,\ \dot{z} = G(e,z)
\end{equation}
with $e$ in $\RR^{(N-1)n}$, $z$ in $\RR^n$ and where
\begin{align}
\label{eq_VecFieldSynchro1}F &= \bbm{F_2^\top & F_3^\top & \dots & F_N^\top}^\top\\
F_i(e,z)&=
 f(z+e_i)-f(z) \label{eq_VecFieldSynchro2}
\\
\nonumber&\qquad
+ g(z+e_i)\bar \phi_i(e,z)-g(x_1)\bar \phi_1(e,z)\ ,
\\[0.3em]
\label{eq_VecFieldSynchro3}G(e,z)&=
 f(z)+g(z)\bar \phi_1(e,z) \ ,
\end{align}
where we have used the notation
\begin{equation}\label{LP22}
\bar \phi_i(e,z) = \phi_i(z,z+e_2,\dots, z+e_N).
\end{equation}
\startmodifVA
Note that we have
\begin{align}
\nonumber |e|^2 &= \sum_{i=2}^N |x_i-x_1|^2\ ,\\
&\leq 2\sum_{i=2}^N |x_i-\bar x|^2 +2 (N-1) |\bar x- x_1|^2\\
\label{eq_EqNorme1}
&\leq 2(N-1)|x|_\DR^2\ ,
\end{align}
and,
\begin{align}
\nonumber
|x|_\DR^2 
&= \min_{z\in \RR^n}\sum_{i=1}^N |z-x_i|^2\\
 &\leq \sum_{i=1}^N |x_1-x_i|^2
= |e|^2.\label{eq_EqNorme2}
%
\end{align}

Hence, if we denote $E(e,z,t)$ the $e$ components of the solution to (\ref{eq_System}), then (\ref{eq_Beta_exp}) implies for all $(e,z)$ in $\RR^{(N-1)n}\times\RR^n$
\begin{multline*}
|E(e,z,t)|\\\leq  \sqrt{(N-1)\left  (1+\frac{N-1}{N^2}\right )} k\exp(-\lambda t)\,  |e|.
\end{multline*}
It follows that the manifold $e=0$ is locally uniformly (in $z$) exponentially stable for (\ref{eq_System}). 
In other words, property \textsf{TULES-NL} (see Section \ref{sec_Trans} in the Appendix) is satisfied.
Employing the assumptions on the bounds on $f$, $g$, $\phi$ and its derivatives, we conclude with \cite[Prop. 1]{AndrieuJayawardhanaPraly_TAC_TransExpStab} that the so-called Property  \textsf{ULMTE}  is
satisfied (see Section \ref{sec_Trans} in the Appendix for the definition). 
Hence there exists a $C^1$ function with matrix valued $P_N:\RR^{n}\rightarrow\RR^{(N-1)n\times(N-1)n}$ and a positive definite matrix $Q_N$ in $\RR^{(N-1)n\times(N-1)n}$ such that for all $z$ in $\RR^n$
\begin{multline}\label{eq_VA3}
\mathfrak{d}_f P_N(z) + P_N(z)\frac{\partial F}{\partial e}(0,z) + \frac{\partial F}{\partial e}(0,z)^\top \leq -Q_N,
\end{multline}
and
\begin{equation}\label{eq_VA}
\underline p_{N} \Id \leq P_N(z) \leq \overline p_N\Id\ .
\end{equation}
For each $z$, let us decompose
$$
P_N(z) = \begin{bmatrix}
S(z) & T(z)\\T(z)^\top & R(z) 
\end{bmatrix},
$$
with $S$ taking value in $\RR^{n\times n}$ and $T$ and $R$ of appropriate dimensions.

Consider the $C^1$ matrix function $P:\RR^{n}\rightarrow\RR^{n\times n}$ defined as follows.
\begin{align*}
P(z) &=  \begin{bmatrix}\Id & -T(z) R(z)^{-1}\end{bmatrix} P_N(z)\begin{bmatrix}\Id \\-R(z)^{-1}T(z)^\top\end{bmatrix},\\
&=S(z) - T(z)R(z)^{-1}T(z)^\top\ .
\end{align*}
We will show that this matrix function $P$ satisfies all assumptions of property \textsf{CMF}.
First of all, we show that $P$ satisfies (\ref{LP10}). 
Pre- and post- multiplying equation (\ref{eq_VA}) by the two matrices $\begin{bmatrix}\Id & -T(z) R(z)^{-1}\end{bmatrix}$ and
$\begin{bmatrix}\Id \\-R(z)^{-1}T(z)^\top\end{bmatrix}$, yields 
\begin{multline*}
\underline p_{N} (\Id +T(z) R(z)^{-2}T(z)^\top ) \leq P(z) \\
\leq \overline p_N(\Id +T(z) R(z)^{-2}T(z)^\top ).
\end{multline*}
On another hand,
$$
|T(z)T(z)^\top|  = \left|\begin{bmatrix}
0 & \Id
\end{bmatrix}P(z)\begin{bmatrix}
\Id &0\\0&0
\end{bmatrix}P(z)\begin{bmatrix}0\\\Id\end{bmatrix}\right|\leq \overline p_N^2.
$$
Moreover,
$$
\underline p_N\Id \leq R(z)\leq  \overline p_N\Id.
$$
Which gives by pre and post multiplying by $R(z)^{-\frac{1}{2}}$
$$
R(z)^{-1} \underline p_N \leq \Id \leq \overline p_N R(z)^{-1}.
$$
Consequently, it yields equation (\ref{LP10}) since we have
$$
\underline p_{N} \Id
\leq P(z) \leq \overline p_N\left(1+ \frac{\overline p_N^2}{\underline p_N^2} \right) \Id.
$$

We now show that (\ref{eq_TensorDerivSynch}) holds.
Note that we have
\begin{multline*}
\mathfrak d_f P(z) = \mathfrak d_f S(z)  - \mathfrak d_f T(z)R(z)^{-1}T(z)^\top \\ 
+ T(z)R(z)^{-1}\mathfrak d_f R(z)R(z)^{-1}T(z)^\top \\ - T(z)R(z)^{-1}\mathfrak d_f T(z)^\top ,
\end{multline*}
which gives
\begin{equation}\label{eq_VA7}
\mathfrak d_f P(z)  = \begin{bmatrix}\Id & -T(z) R(z)^{-1}\end{bmatrix}\mathfrak d_f  P_N(z)\begin{bmatrix}\Id \\-R(z)^{-1}T(z)^\top\end{bmatrix}.
\end{equation}

Note now that
\begin{multline}
\label{eq_VA2}
\begin{bmatrix}\Id & -T(z) R(z)^{-1}\end{bmatrix} P_N(\dx)
\frac{\partial F}{\partial \de}(0,\dx)\begin{bmatrix}\Id \\-R(z)^{-1}T(z)^\top\end{bmatrix}
\\
\begin{aligned}
&= \begin{bmatrix} P(z) & 0 \end{bmatrix} 
\frac{\partial F}{\partial \de}(0,\dx)\begin{bmatrix}\Id \\-R(z)^{-1}T(z)^\top\end{bmatrix}
\\&= P(z)\frac{\partial F_2}{\partial \de}(0,\dx)\begin{bmatrix}\Id \\-R(z)^{-1}T(z)^\top\end{bmatrix}\\
\end{aligned}
\end{multline}

On another hand, by the definition of $\bar \phi$ in (\ref{LP22}) and the second point of Definition \ref{Def1}, it follows that  $\bar \phi_i(0,z)=0$. 
This implies that for every $i=2,\ldots,N$, 
\begin{multline}\label{eq_dFidi}
\frac{\partial F_i}{\partial \de_i}(0,\dx)\;=\;
\frac{\partial f}{\partial \dx}(\dx)
\\+
g(\dx) \left[\frac{\partial \bar\phi_i}{\partial \de_i}(0,\dx)
-
\frac{\partial \bar\phi_1}{\partial \de_{i}}(0,\dx)
\right]\ 
\end{multline}
and for all $j\neq i$, 
\begin{equation}\label{eq_dFidj}
\frac{\partial F_i}{\partial \de_{j}}(0,\dx)\;=\;
g(\dx) \left[\frac{\partial \bar\phi_i}{\partial \de_{j}}(0,\dx)
-
\frac{\partial \bar\phi_1}{\partial \de_{j}}(0,\dx)
\right]\ .
\end{equation}

Consequently,
\begin{equation}\label{eq_VA4}
\frac{\partial F_2}{\partial \de}(0,\dx)\begin{bmatrix}\Id \\-R(z)^{-1}T(z)^\top\end{bmatrix} = \frac{\partial f}{\partial \dx}(\dx) + g(\dx)h(\dx)
\end{equation}
where
\begin{multline}\label{eq_VA5}
h(\dx) = \begin{bmatrix}
\frac{\partial \bar\phi_2}{\partial \de_{2}}(\dx)
-
\frac{\partial \bar\phi_1}{\partial \de_{2}}(\dx)
&
\dots
&
\frac{\partial \bar\phi_2}{\partial \de_{N}}(\dx)
-
\frac{\partial \bar\phi_1}{\partial \de_{N}}(\dx)
\end{bmatrix}\\\times\begin{bmatrix}\Id \\-R(z)^{-1}T(z)^\top\end{bmatrix} .
\end{multline}
%
Consequently,  pre- and post- multiplying equation (\ref{eq_VA3}) by the two matrices $\begin{bmatrix}\Id & -T(z) R(z)^{-1}\end{bmatrix}$ and
$\begin{bmatrix}\Id \\-R(z)^{-1}T(z)^\top\end{bmatrix}$, and employing  equations (\ref{eq_VA7}), (\ref{eq_VA2}),  (\ref{eq_VA4}) and (\ref{eq_VA5}) yield a positive definite matrix $Q$ such that
\begin{multline}\label{eq_VA6}
\mathfrak d_fP(\dx) + P(\dx)\left[\frac{\partial f}{\partial x}(\dx) + g(\dx)h(\dx)\right]\\
+\left[\frac{\partial f}{\partial x}(\dx) + g(\dx)h(\dx)\right]^\top  P(\dx) \leq - Q.
\end{multline}
From this equation,  (\ref{eq_TensorDerivSynch}) is satisfied and Property \textsf{CMF} holds.

Moreover (\ref{eq_VA6}) implies that property \textsf{ULMTE} introduced in  \cite{AndrieuJayawardhanaPraly_TAC_TransExpStab}  (see Appendix \ref{sec_Trans}) is satisfied for the system
$$
\dot z = \bar F(e,z)\ ,\ \dot z = \bar G(e,z),
$$
where
$$
\bar F(e,z) = f(e+z) - f(z) + g(\dx)h(\dx)e\ ,\ \bar G(e,z) = f(z).
$$
Hence, employing Proposition \ref{Prop_ExistTensor} in the appendix, one can conclude that property \textsf{IS} is satisfied with the control $\tilde u = h(z)\tilde z$.
\end{proof}
\fi

In the following section, we discuss the possibility to design an exponential synchronizing control law based on these necessary conditions.
\section{Sufficient condition}
\label{Sec_subSuffCond}
\subsection{Sufficient conditions for local exponential synchronization}
The interest of the Property \textsf{CMF} given in Subsection \ref{subsection_necessary} is to use the symmetric covariant tensor $\PR$ in the design of a local synchronizing control law. 
Indeed, following one of the main results in \cite{AndrieuJayawardhanaPraly_TAC_TransExpStab}, we get the following sufficient condition for the  solvability of (local) uniform exponential synchronization problem.
The first assumption is that, up to a scaling factor, the control vector field $g$ is a gradient field with $P$ as a Riemannian metric (see also \cite{Fornietal_CDC_2013} for similar integrability assumption). The second one is related to the \textsf{CMF} property.
\begin{theorem}[Local sufficient condition]
\label{Theo_SufCond}
Assume that $g$ is bounded  and that $f$ and $g$ have bounded first and second derivatives.
Assume that there exists a $C^2$ function $\PR :\RR^n\rightarrow\RR^{n\times n}$
which values are symmetric positive definite matrices and
with a bounded derivative that satisfies the following two conditions.  
\begin{list}{}{%
\parskip 0pt plus 0pt minus 0pt%
\topsep 0pt plus 0pt minus 0pt
\parsep 0pt plus 0pt minus 0pt%
\partopsep 0pt plus 0pt minus 0pt%
\itemsep 0pt plus 0pt minus 0pt
\settowidth{\labelwidth}{1.}%
\setlength{\labelsep}{0.5em}%
\setlength{\leftmargin}{\labelwidth}%
\addtolength{\leftmargin}{\labelsep}%
}
\item[1.]
There exist a $C^2$ function $U:\RR^n\to \RR$ which
has bounded first and second derivatives,
and a $C^1$ function $\alpha :\RR^n\to \RR^p$ which has bounded first and second derivatives such
that
\begin{equation}
\label{LP24}
\frac{\partial U}{\partial \dx}(\dx)^\top = P(\dx)g(\dx)\alpha(\dx)
\; ,
\end{equation}
holds for all $\dx$ in $\RR^n$; and
\item[2.]
There exist a symmetric positive definite matrix $Q$ and positive constants 
$\underline{p}$, $\overline{p}$  and $\rho >0$
such that (\ref{LP10}) holds
and
\begin{equation}
\label{LP23}
\hskip -1em
 L_f\PR(\dx)
- \rho 
\frac{\partial U}{\partial \dx}(\dx)^\top \frac{\partial U}{\partial \dx}(\dx)
\leq  - Q
\  ,
\end{equation}
 hold for all $\dx$ in $\RR^n$. 
\end{list}
Then, given a connected graph $\mathcal G$ with associated Laplacian matrix $L=(L_{ij})$, there exists a  constant  $\underline{\kell}$
 such that the control law $u=\phi(x)$ with $\phi=\bbm{\phi_1^\top & \dots & \phi_N^\top}^\top$  given by 
\begin{equation}\label{eq_ControlLaw2}
\phi_i(x)=-\kell\alpha (x_i)\,  \sum_{j=1}^N L_{ij} U(x_j)
\end{equation}
with $\kell\geq \underline{\kell}$  solves the local uniform exponential synchronization of (\ref{eq_Syst}). 
\end{theorem}
\if\arxiv0
The proof of this result can be found in \cite{AndrieuJayawardhanaTarbouriech_CDC_Synchro} or in \cite{AndrieuJayawardhanaTarbouriech_arXiv_16_Synchro}.
\fi
\begin{remark}
Assumption (\ref{LP23}) is stronger than the necessary condition \textsf{CMF}. Note however, that employing some variation on Finsler Lemma (see \cite{AndrieuJayawardhanaPraly_TAC_TransExpStab} for instance) it can be shown that these assumptions are equivalent when $x$ remains in a compact set.
\end{remark}
\begin{remark}Note that 
 for all $x=\mathds{1}_N\otimes \bx=(\bx,\dots,\bx)$ in $\DR$ and for all $(i,j)$ with $i\neq j$
\begin{equation}
\frac{\partial \phi_i}{\partial x_j}(x) = -\ell \alpha(\bx) L_{ij} \frac{\partial U}{\partial \bx}(\bx) .
\end{equation}
Hence, for all $x=\mathds{1}_N\otimes \bx$ in $\DR$, we get 
\begin{equation}\label{eq_dphidw}
\frac{\partial \phi}{\partial x}(x) = -\ell L \otimes  \alpha(\bx) \frac{\partial U}{\partial \bx}(\bx)  \ .
\end{equation}
\end{remark}
\if\arxiv1
\begin{proof}
First of all, note that the control law $\phi$ satisfies the condition 1) and 2) in Definition \ref{Def1}.
Indeed, for all $x$ and all $(i,j)$ with $i\neq j$
$$
\frac{\partial \phi_i}{\partial x_j}(x) = -\ell \alpha(x_i) L_{ij} \frac{\partial U}{\partial \bx}(x_j) \ .
$$
If $(i,j)\notin \mathcal E$, it yields $L_{ij}=0$ and consequently $\frac{\partial \phi_i}{\partial x_j}(x)=0$.

Moreover, when $x$ is in $\DR$, i.e., $x=\mathds{1}_N\otimes \bx=(\bx,\dots,\bx)$ for all $i$
$$
\phi_i(x)=-\kell\alpha (z)\,  \left(\sum_{j=1}^N L_{ij}\right) U(z) =0.
$$

It remains to show that condition 3) of Definition \ref{Def1} holds. More precisely, we need to prove that the manifold $\mathcal D$  is locally exponentially stable along the solution of  the closed-loop system. 

As in the proof of Theorem \ref{theo_Necessary}, let us denote $e=(e_{2},\dots, e_{N})$ with $e_{i}=x_1-x_i$ and $z=x_1$.  
Note that the closed-loop system may be rewritten as in (\ref{eq_System})
with the vector fields $F$ and $G$ as defined in (\ref{eq_VecFieldSynchro1})--(\ref{eq_VecFieldSynchro3}) with $\phi$ as the control law.

The rest of the proof is to apply \cite[Proposition 3]{AndrieuJayawardhanaPraly_TAC_TransExpStab}. 
For this purpose, we need to show that for closed-loop system (\ref{eq_VecFieldSynchro1})--(\ref{eq_VecFieldSynchro3}) the property \textsf{ULMTE} introduced in  \cite{AndrieuJayawardhanaPraly_TAC_TransExpStab} and given in Section \ref{sec_Trans} is satisfied.

By the assumption on the graph being connected and together with Lemma \ref{lem_Laplacian}, we have that the matrix
$
A = -(\Lc - \mathds{1}_{N-1}\Lb)
$
is Hurwitz. 
Let $S$ in $\RR^{(N-1)\times(N-1)}$ be a symmetric positive definite matrix solution to the Lyapunov equation
\begin{equation}\label{eq_Lyap}
SA  + A^\top S \leq -\nu S
\end{equation}
where $\nu$ is a positive real number.

Consider the $C^1$ function $P_N:\RR^n\rightarrow\RR^{(N-1)n\times(N-1)n}$ defined as
$$
P_N(\dx) =S\otimes P(\dx).
$$
Our aim is to show that the closed loop system satisfies  property \textsf{ULMTE} given in Section \ref{sec_Trans}.
First of all, note that $S$ being symmetric positive definite, with (\ref{LP10}), it yields the existence of positive real numbers $\underline{p}_N$, $\overline p_N$ such that
$$
\underline p_N \Id_{N-1} \leq P_N(z) \leq \overline p_N \Id_{N-1}\ .
$$
Hence, equation (\ref{LP7}) is satisfied.

Note that we have $G(0,\dx) = f(\dx)$. Moreover we have
$$
\der_{G(0,\dx)} P_N(\dx)  = S \otimes \der_{f}P(\dx).
$$

Note that with properties (\ref{eq_dFidi}), (\ref{eq_dFidj}) and (\ref{eq_dphidw}), it follows that 
\begin{multline}
\frac{\partial F}{\partial \de}(0,\dx) = \Id_{N-1}\otimes \frac{\partial f}{\partial \dx}(\dx) \\+ \ell
A \otimes\left(\alpha(\dx)g(\dx)\frac{\partial U}{\partial \dx}(\dx)\right).
\end{multline}
Hence, 
\\[0.5em]$\displaystyle
\der_{G(0,\dx)} P_N(\dx) +
P_N (\dx ) \frac{\partial F}{\partial \de}(0,\dx)
+ \frac{\partial F}{\partial \de}(0,\dx)^\top P_N(\dx) 
$\\[0.5em]\null\ $\displaystyle
=S\otimes \left(\der_{f}P(\dx)+P(\dx)\frac{\partial f}{\partial \dx}(\dx)+\frac{\partial f}{\partial \dx}(\dx)^\top P(\dx)\right)
$\hfill\null\\[0.5em]\null\hfill$\displaystyle
+\ell (SA+A^\top S)\otimes \left (\frac{\partial U}{\partial \dx}(\dx)^\top\frac{\partial U}{\partial \dx}(\dx)\right ).
$\\[0.5em]
With (\ref{eq_Lyap}) and (\ref{LP23}) this implies that 
\begin{multline*}
\der_{f} P_N(\dx) +
P_N (\dx ) \frac{\partial F}{\partial \de}(0,\dx)
+ \frac{\partial F}{\partial \de}(0,\dx)^\top P_N(\dx)\\ \leq 
S\otimes \left (-Q +(\rho-\ell \nu) \frac{\partial U}{\partial \dx}(\dx)^\top\frac{\partial U}{\partial \dx}(\dx)\right ).
\end{multline*}
Hence, by choosing $\underline{\kell}\geq \frac{\rho}{\nu}$,  inequality (\ref{eq_TensorDerivative}) holds and consequently Property \textsf{ULMTE} holds.
The last part of the proof is to make sure that the vector field $F$ has bounded first and second derivatives and that the vector field $G$ has bounded first derivative. 
Note that by employing the bounds on the functions $P$, $f$, $g$, $\alpha$ and their derivatives, the result immediately follows from Proposition \ref{Prop_Lyap} in Section \ref{sec_Trans}.  
Indeed, this implies that Property \textsf{TULES-NL} holds and consequently,  $e=0$ is (locally) exponentially stable manifold  
for system (\ref{eq_VecFieldSynchro1})--(\ref{eq_VecFieldSynchro3}) in closed loop with the control (\ref{eq_ControlLaw2}).
With inequalities (\ref{eq_EqNorme1}) and (\ref{eq_EqNorme2}), it implies that inequality (\ref{eq_Beta_exp}) holds for $r$ sufficiently small. 
\end{proof} 
\else

\fi

\subsection{Sufficient conditions for global exponential synchronization}

Note that in \cite{AndrieuJayawardhanaPraly_TAC_TransExpStab} with an extra assumption related to the metric (the level sets of $U$ are totally geodesic sets with respect to the Riemannian metric obtained from $\PR$), it is shown that global synchronization may be achieved when considering only two agents which are connected. It is still an open question to know if global synchronization may be  achieved  in the general nonlinear context with more than two agents.
However in the particular case in which the matrix $\PR(z)$ and the vector field $g$ are constant, then global synchronization may be achieved as this is shown in the following theorem.

\begin{theorem}[Global sufficient condition]
\label{Theo_SufCondGlob2}
Assume  that $g(z)=G$ and there exists a symmetric positive definite matrix $\PR$ in $\RR^{n\times n}$,
 a symmetric positive definite matrix $Q$ and $\rho >0$ such that
\begin{equation}
\label{LP23Eucl}
\PR\frac{\partial f}{\partial \dx}(\dx) + \frac{\partial f}{\partial \dx}(\dx)^\top \PR
- \rho  \PR GG^\top \PR \leq  -Q
\  .
\end{equation}
Assume moreover that the graph is connected with Laplacian matrix $L$.
Then there exist  constants  $\underline{\kell}$ and positive real numbers $c_1, \dots, c_N$
 such that the control law $u=\phi(x)$ with $\phi=\bbm{\phi_1^\top & \dots & \phi_N^\top}^\top$  given by 
\begin{equation}\label{eq_ControlLaw}
\phi_i(x)=-\kell \, c_i \sum_{j=1}^N  L_{ij} G^\top P x_j
\end{equation}
with $\kell\geq \underline{\kell}$,  solves the global uniform exponential synchronization for (\ref{eq_Syst}). 
\end{theorem}
\begin{proof}
Let $c_j=1$ for $j=2,\dots,N$. Hence only $c_1$ is different from $1$ and remains to be selected.
\if\arxiv1
As in the proof of Theorem \ref{theo_Necessary}, 
let us denote $e=(e_{2},\dots, e_{N})$ with $e_{i}=x_1-x_i$ and $z=x_1$.  
\else
Let us denote $e=(e_{2},\dots, e_{N})$ with $e_{i}=x_1-x_i$ and $z=x_1$.  
\fi
Note that for $i=2,\dots, N$, we have along the solution of the system (\ref{eq_Syst}) with $u$ defined in (\ref{eq_ControlLaw}),
\begin{multline*}
\dot e_i = f(z)-\kell \,   c_1 \sum_{j=1}^N  L_{1j} GG^\top P x_j\\
 -f(z+e_i)+\kell \,    \sum_{j=1}^N  L_{ij} GG^\top P x_j.
\end{multline*}
Note that $L$ being a Laplacian, we have for all $i$ in $[1,N]$ the equality $\sum_{j=1}^N  L_{ij}=0$.
Consequently, we can add  the term $\ell c_1\sum_{j=1}^N  L_{1j} GG^\top P x_1$ and substract the term $\ell\sum_{j=1}^N  L_{ij} GG^\top P x_1$ in the  preceding equation above  so that for $i=2,\dots,N$
\begin{align*}
\dot e_i & = f(z)-\kell \,   c_1 \sum_{j=1}^N  L_{1j} GG^\top P (x_j-x_1) \\ & \qquad 
 -f(z+e_i)+\kell \,    \sum_{j=1}^N  L_{ij} GG^\top P (x_j-x_1), \\
& = f(z)-f(z+e_i)-\kell \,    \sum_{j=2}^N  \left(L_{ij}-c_1L_{1j}\right) GG^\top P e_j.
\end{align*}
 One can check that these equations can be written compactly as  
$$
\dot e = \left[\int_0^1\Delta(z,e,s) ds +\kell \,  (A(c_1)\otimes GG^\top P) \right]e,
$$
with $A(c_1)$ is matrix in $\RR^{(N-1)\times (N-1)}$, which depends on the parameter $c_1$ and  is obtained from the Laplacian as~:
$$
A(c_1) = -\left[L_{2:N,2:N} -c_1 L_{1,2:N}\mathds{1}_{N-1} \right],
$$
where
 $L= \begin{bmatrix}
\La & \Lb\\
\Lb^\top & \Lc
\end{bmatrix}$
and $\Delta$ is the $(N-1)n\times n$ matrix valued function defined as
$$
\Delta(z,e,s)=\texttt{Diag}\left\{
\frac{\partial f}{\partial z}(z-se_2),\dots, \frac{\partial f}{\partial z}(z- se_N)\right\}
\ .
$$
The following Lemma shows that by selecting $c_1$ sufficiently small the matrix $A$ satisfies the following property.
Its proof is given in the Appendix.

\begin{lemma}\label{lem_Laplacian2}
If the communication graph is connected then there exist
sufficiently small $c_1$   and $\mu>0$ such that
$$
A(c_1) + A(c_1)^\top \leq -\mu I
$$
\end{lemma}
\vspace{0.2cm}

With this lemma in hand, we consider now the candidate Lyapunov function  
$
V(e) = e^\top P_N e\ ,
$
where $P_N = (I_{N-1} \otimes P)$.
Note that along the solution, the time derivative of this function satisfies~:
$$
\dot{\overparen{V(e)}} = 2e^\top P_N\left[\int_0^1\Delta(z,e,s) ds +\kell \,  (A(c_1)\otimes GG^\top P) \right]e\ .
$$
Note that we have
\begin{multline*}
P_N\Delta(z,e,s) \\= \texttt{Diag}\left\{
P\frac{\partial f}{\partial z}(z-se_2), \dots,P\frac{\partial f}{\partial z}(z- se_N)
\right\}\ ,
\end{multline*}
and
\begin{align*}
2e^\top(I_{N-1} \otimes P)&(A(c_1)\otimes GG^\top P) e \\
&=
e^\top( [A(c_1)+A(c_1)^\top]\otimes P GG^\top P) e\ \\
&\leq 
-e^\top( \mu I_{N-1}\otimes P GG^\top P )e\ .
\end{align*}
Hence, we get
$
\dot{\overparen{V(e)}} \leq \int_0^s e^\top M(e,z,s) e \, ds
$,
where $M$ is the  $(N-1)n\times (N-1)n$ matrix defined as
$$
M(e,z,s)= \texttt{Diag}\left\{
M_2(e,z,s), \dots, M_{N}(e,z,s)
\right\} \ ,
$$
with, for $i=2, \dots, N$
\begin{multline*}
M_i(e,z,s) = 
P\frac{\partial f}{\partial z}(z-se_i)+\frac{\partial f}{\partial z}(z-se_i)^\top P \\-2\kell \mu P GG^\top P.
\end{multline*}
Note that by taking $\kell$ sufficiently large, with (\ref{LP23Eucl}) this yields  $M_i(e,z,s)\leq -Q$.  This immediately implies that  
$
\dot{\overparen{V(e)}} \leq -e^\top\left( I_{N-1}\otimes Q\right) e
$.
This ensures exponential convergence of $e$ to zero on the time of existence of the solution.
\if\arxiv1
With (\ref{eq_EqNorme1}) and (\ref{eq_EqNorme2}), this yields global exponential synchronization of the closed-loop system.
\else
\startmodifVA
Let $\bar x = \texttt{argmin}_{z\in\RR^n} \sum_{i=1}^N |z - x_i|^2$.
Note that we have
\begin{align}
\nonumber |e|^2 
&\leq 2\sum_{i=2}^N |x_i-\bar x|^2 +2 (N-1) |\bar x- x_1|^2\\
&\leq 2(N-1)|x|_\DR^2\ ,
\end{align}
\begin{align}
\nonumber
|x|_\DR^2 
&= \min_{z\in \RR^n}\sum_{i=1}^N |z-x_i|^2\\
 &\leq \sum_{i=1}^N |x_1-x_i|^2
= |e|^2.\label{eq_EqNorme2}
%
\end{align}
This yields global exponential synchronization of the closed-loop system.
\stopmodif
\fi
\end{proof}

In the following section, we show that the property CMF required to design a distributed synchronizing control law can be obtained for a large class of nonlinear systems. This is done via backstepping design.

\section{Construction of an admissible tensor via backstepping}
\label{Sec_BackStep}
\subsection{Adding derivative (or backstepping)}

As proposed in Theorem \ref{Theo_SufCond}, a distributed synchronizing control law can be designed using a symmetric covariant tensor field of order 2, which satisfies (\ref{eq_TensorDerivSynch}).
Given a general nonlinear system, the construction of such a matrix function $P$ may be a hard task. 
In \cite{sanfelice2015solution}, a construction of the function $P$ for observer based on the integration of a Riccati equation is introduced. 
Similar approach could be used in our synchronization problem.
Note however that in our context an integrability condition (i.e. equation (\ref{LP24})) has to be satisfied by the function $P$.
This constraint may be difficult to address when considering a Riccati equation approach.

In the following we present a constructive design of such a matrix $P$ that resembles the backstepping method. 
This approach can be related to \cite{ZamanietAl_SCL_2013backstepping,ZamaniTabuada_TAC_11BacksteppingIncrStab} in which a metric is also constructed iteratively. 
We note that one of the difficulty we have here is that we need to propagate the integrability property given in equation (\ref{LP24}). 

For outlining the backstepping steps for designing $P$, we consider the case in which the vector fields $(f,g)$ can be decomposed as follows
\begin{equation*}
f(z) = \begin{bmatrix}
\fa (\bxa) + g_a(\bxa)\bxb\\
\fb(\bxa,\bxb)
\end{bmatrix}\ ,
\end{equation*}
and,
\begin{equation*}
g(\bx) = \begin{bmatrix}
0\\ \gb(\bx) 
\end{bmatrix}\ ,\ 0<\underline{g}_b\leq \gb(\bx) \leq \overline{g}_b
\end{equation*}
with $\bx=\bbm{\bxa^\top & \bxb}^\top$, $\bxa$ in  $\RR^{\dna}$ and $\bxb$ in $\RR$. In other words, 
\begin{equation}\label{backstepping_sys_eq}
\dot{\bx}_a = \fa (\bxa) + \ga (\bxa)\bxb, \ \dot{\bx}_b = f_b(\bx) + g_b(\bx) u.
\end{equation}
Let $\Ca $ be a compact subset of $\RR^{\dna}$. As in the standard backstepping approach, we make the following assumptions on the $\bxa$-subsystem where $\bxb$ is treated as a control input to this subsystem. 
\begin{assumption}[$\bxa$-Synchronizability]
\label{Ass_x1Sync}
 Assume that there exists a $C^\infty$ function $\PRa  :\RR^{\dna}\rightarrow\RR^{\dna\times \dna}$
that satisfies the following conditions. 
\begin{list}{}{%
\parskip 0pt plus 0pt minus 0pt%
\topsep 0pt plus 0pt minus 0pt
\parsep 0pt plus 0pt minus 0pt%
\partopsep 0pt plus 0pt minus 0pt%
\itemsep 0pt plus 0pt minus 0pt
\settowidth{\labelwidth}{1.}%
\setlength{\labelsep}{0.5em}%
\setlength{\leftmargin}{\labelwidth}%
\addtolength{\leftmargin}{\labelsep}%
}
\item[1.]
There exist a $C^\infty$ function $\Ua :\RR^{\dna}\to \RR$ 
and a $C^\infty$ function $\alphaa:\RR^{\dna}\to \RR$ such
that
\begin{equation}
\label{LP241}
\frac{\partial \Ua }{\partial \bxa}(\bxa)^\top = \alphaa (\bxa)\PRa (\bxa)\ga (\bxa)
\end{equation}
holds for all $\bxa$ in $\Ca$;
\item[2.]
There exist a symmetric positive definite matrix $\Qa$ and  positive constants 
$\underline{p}_a$, $\overline{p}_a$  and $\rhoa  >0$
such that 
\begin{equation}
\underline{p}_a\Id_{n_a}\leq \PRa(\bxa)\leq \overline{p}_a\Id_{n_a}\ ,\ \forall \bxa\in\RR^{n_a},
\end{equation} 
holds and
\begin{equation}
\label{LP231}
L_{\fa} \PRa  (\bxa)  - \rhoa  
\frac{\partial \Ua }{\partial \bxa}(\bxa)^\top
\frac{\partial \Ua }{\partial \bxa}(\bxa)\leq  - \Qa  
\  ,
\end{equation}
%
%
 holds for all $\bxa$ in $\Ca$. 
\end{list}
\end{assumption}

As a comparison to the standard backstepping method for stabilizing nonlinear systems in the strict-feedback form, the $\bxa$-synchronizability conditions above are akin to the stabilizability condition of the upper subsystem via a control Lyapunov function. However, for the synchronizer design as in the present context, we need an additional assumption to allow the recursive backstepping computation of the tensor $P$.   
Roughly speaking, we need the existence of a mapping $\qa$ such that the metric $\PR_a $ becomes invariant along the vector field $\frac{\ga}{\qa}$.
In other words,  $\frac{\ga}{\qa}$ is a Killing vector field.

\begin{assumption}
\label{Ass_Pinv}
There exists a non-vanishing smooth function $\qa :\RR^{\dna}\rightarrow\RR$ such that the metric obtained from $\PR_a $ on $\Ca $ is invariant along $\frac{\ga (\bxa)}{\qa(\bxa)}$. 
In other words, for all $\bxa$ in $\Ca$
\begin{equation}\label{eq_AssPinv}
 L_{\frac{\ga (\bxa) }{\qa (\bxa) }}\PRa(\bxa)  =0\ .
\end{equation} 
\end{assumption}
\vspace{0.1cm}

Similar assumption can be  found in \cite{Fornietal_CDC_2013} in the characterization of differential passivity.

Based on the Assumptions \ref{Ass_x1Sync} and \ref{Ass_Pinv}, we have the following theorem on the backstepping method for constructing a symmetric covariant tensor field $\PRb$ of the complete system (\ref{backstepping_sys_eq}). 

\begin{theorem}\label{Theo_Backstepping}
Assume that the $\bxa$-subsystem satisfies Assumption \ref{Ass_x1Sync} and Assumption  \ref{Ass_Pinv}  in the compact set $\Ca$
with a $\dna\times \dna$ symmetric covariant tensor field $\PR_a$ of order two  and a non-vanishing smooth mapping $\qa:\RR^{\dna}\rightarrow\RR$. 
Then for all positive real number $M_b$, the system (\ref{backstepping_sys_eq}) with the state variables $\bx=(\bxa,\bxb)\in\RR^{\dna+1}$ satisfies the Assumption \ref{Ass_x1Sync} in the compact set $\Ca \times [-M_b,M_b]\subset \RR^{\dna+1}$ with the symmetric covariant tensor field $\PRb$ be given by
\begin{equation*}
\PRb (\bx) = 
\begin{bmatrix}
\PR_a (\bxa) + \Sa (\bx) \Sa (\bx)^\top
&  \Sa (\bx)\qa(\bxa)	\\
\Sa (\bx)^\top \qa (\bxa) & \qa (\bxa)^2
\end{bmatrix}
\end{equation*}
where
$
\Sa (\bx) = \frac{\partial \qa }{\partial \bxa}(\bxa)^\top \bxb+\eta \alphaa (\bxa)\PR_a (\bxa) \ga (\bxa)
$
and $\eta$ is a positive real number. Moreover, there exists a non-vanishing mapping $\qb:\RR^{\dna+1}\rightarrow\RR$ such that $\PRb$ is invariant along $\frac{g}{\qb}$. 
In other words, Assumptions \ref{Ass_x1Sync} and \ref{Ass_Pinv} hold for the complete system (\ref{backstepping_sys_eq}). 
\end{theorem}
\vspace{0.1cm}
\begin{remark}
Note that with this theorem, since we propagate the required property we are able to obtain a synchronizing control law for any triangular nonlinear system.
\end{remark}
\begin{proof}
Let $M_b$ be a positive real number and let $\Cb=\Ca \times [-M_b,M_b]$.
Let $\Ub :\RR^{\dna+1}\rightarrow\RR$ be the function defined by
$$
\Ub (\bxa,\bxb) =  \eta \Ua (\bxa) + \qa (\bxa)\bxb\ .
$$
where $\eta$ is a positive real number that will be selected later on.
It follows from (\ref{LP241}) that for all $(\bxa,\bxb)\in \Cb$, we have 
\begin{align*}
\frac{\partial \Ub}{\partial \bx}(\bx) ^\top 
&= \begin{bmatrix}
\eta \frac{\partial \Ua }{\partial \bxa}(\bxa)^\top+\frac{\partial \qa }{\partial \bxa}(\bxa)\bxb\\\qa (\bxa)
\end{bmatrix}\\
&=\frac{1}{\qa (\bxa)}\PRb (\bx) \begin{bmatrix}
0\\
 1
\end{bmatrix}\\&=\alphab(\bx) \PRb (\bx) g(\bx) 
\end{align*}
with $\alphab(\bx) =\frac{1}{\qa (\bxa)\gb(\bx)}$. Hence, the first condition in  Assumption \ref{Ass_x1Sync} is satisfied. 

Consider $\bx$ in $\Cb$ and let $v=\bbm{\va^\top & \vb}^\top$ in $\RR^{n_a+1}$ be such that
\begin{equation}\label{eq_LgVzero}
v^\top\PRb(\bx) g(\bx) =0\ .
\end{equation}
Note that this implies that 
\begin{equation}\label{eq_LgVzero2}
\vb  = -\va ^\top \frac{\Sa (\bx)}{\qa (\bxa)}.  
\end{equation}
In the following, we compute the expression~:
$$v^\top L_f \PRb(\bx)v=v^\top \der_f\PRb  (\bx) v+
2v^\top \PRb  (\bx) \frac{\partial f}{\partial \bx}(\bx) v\ .$$
For the first term, 
we have
\begin{multline*}
v^\top \der_f\PRb  (\bx) v =\va ^\top \der_{\fa }\PRa  (\bxa)\va  + \bxb \va ^\top \der_{\ga }\PRa  (\bxa)\va\\
 + \va ^\top \der_{f} \Sa (\bx) \Sa (\bx)^\top \va  +2 \va ^\top \der_{f } \Sa (\bx)\qa (\bxa)\vb  \\+ \der_{\fa +\ga \bxb} \qa  (\bxa)^2 \vb^2 
\end{multline*}
With (\ref{eq_LgVzero2}), it yields
\begin{multline*}
 \va ^\top \der_{f} \Sa (\bx) \Sa (\bx)^\top \va  +2 \va ^\top \der_{f } \Sa (\bx)\qa (\bxa)\vb  \\+ \der_{\fa +\ga \bxb} \qa  (\bxa)^2 \vb^2 =0
\end{multline*}
Hence
$$
v^\top \der_f\PRb  (\bx) v =\va ^\top \der_{\fa }\PRa  (\bxa)\va  + \bxb \va ^\top \der_{\ga }\PRa  (\bxa)\va\ .
$$
On the other hand, for the second term we have
$$
\PRb (\bx)  = \begin{bmatrix}
\PRa (\bxa) &0\\0&0
\end{bmatrix} + \frac{\PRb (\bx) g(\bx) g(\bx) ^\top\PRb (\bx) }{(\qa(\bxa)\gb(\bx) )^2}
$$
Hence, with (\ref{eq_LgVzero}), it yields
\\[0.5em]$\displaystyle
v^\top \PRb  (\bx) \frac{\partial f}{\partial \bx}(\bx) v=
\begin{bmatrix}
\va ^\top &-\va ^\top\frac{\Sa (\bx)}{\qa (\bxa)}
\end{bmatrix} 
\PR(\bx) 
$\\[0.5em]\null\hfill$\displaystyle
\begin{bmatrix}
\frac{\partial \fa }{\partial \bxa}(\bxa)+\frac{\partial \ga }{\partial \bxa}(\bxa)\bxb  & \ga (\bxa)\\\frac{\partial \fb}{\partial \bxa}(\bxa,\bxb)&\frac{\partial \fb}{\partial \bxb}(\bxa,\bxb)
\end{bmatrix}\begin{bmatrix}
\va  \\-\frac{\Sa (\bx)^\top}{\qa (\bxa)}\va 
\end{bmatrix}
$\\[0.5em]$\displaystyle \ \ \ \ 
= \va ^\top \PR_a (\bxa)\frac{\partial \fa }{\partial \bxa}(\bxa) \va +\bxb \va ^\top \PRa (\bxa)\frac{\partial \ga }{\partial \bxa}(\bxa) \va  
 $\\[0.5em]$\displaystyle \hphantom{abcabc} - \frac{\eta}{\alphaa (\bxa)\qa(\bxa)} \left|\frac{\partial \Ua }{\partial \bxa}(\bxa)
\va \right|^2$\\[0.5em]
$\displaystyle \hphantom{abcabc} -\frac{\bxb}{\qa(\bxa)}\va ^\top \PRa (\bxa)g(\bxa)\frac{\partial \qa }{\partial \bxa}(\bxa)$\\[0.5em]
Hence, we get
\begin{multline*}
v^\top L_f\PRb (\bx) v= \va ^\top L_{\fa }\PRa (\bxa)\va \\  - \frac{2\eta}{\alphaa (\bxa)\qa(\bxa)} \left|
\frac{\partial \Ua }{\partial \bxa}(\bxa)
\va \right|^2 \\+ \bxb \va ^\top \left [  \der_{\ga }\PRa (\bxa)+\PRa (\bxa)\frac{\partial \ga }{\partial \bxa}(\bxa) \right .\\\left .-2\bxb \va ^\top \PR_a (\bxa)\frac{g(\bxa)}{\qa (\bxa)}\frac{\partial \qa }{\partial \bxa}(\bxa)\right ]\va \ .
\end{multline*}
Let $\eta$ be a positive real number such that
$$
\rhoa  \leq \frac{2\eta}{\alphaa (\bxa)\qa (\bxa)}\ ,\ \forall \bxa\in\Ca \ .
$$
Using (\ref{LP231}) in Assumption \ref{Ass_x1Sync}  and (\ref{eq_AssPinv}) in Assumption \ref{Ass_Pinv}, it follows that for all $\bx$ in $\Cb$ and all $v$ in $\RR^{\dna+1}$ 
\begin{multline*}
v^\top \PRb (\bxa)g(\bx)  =0\\
\Rightarrow v^\top \der_f\PRb  (\bx) v +2v^\top \PRb  (\bx) \frac{\partial f}{\partial x}(\bx) v \leq-v^\top \Qa  v.
\end{multline*}
Employing Finsler theorem and the fact that $\Cb$ is a compact set,
it is possible to show that this implies the existence of a positive real number $\rhob$ such that for all $\bx$ in $\CR_b$
\begin{equation}
L_f\PR(\bx)- \rhob \frac{\partial \Ub }{\partial \bx}(\bx)^\top
\frac{\partial \Ub }{\partial \bx}(\bx) 
 \leq  - \Qb \ .
\end{equation}
where $\Qb$ is a symmetric positive definite matrix.

To finish the proof it remains to show that the metric is invariant along $g$ with an appropriate control law.
Note that if 
$
\qb(\bx) = \qa (\bxa)\gb(\bx)
$
then it follows that this function is also non-vanishing. 
Moreover, we have
\begin{multline*}
L_{\frac{g}{\qb}}\PRb (\bx)= \der_{ \frac{g}{\qb}}\PRb (\bx)  - \frac{\PR(\bx) }{\qa (\bxa)^2}
\begin{bmatrix}
0&0\\
\frac{\partial \qa }{\partial \bxa}(\bxa)&0
\end{bmatrix}
\\-
\begin{bmatrix}
0&\frac{\partial \qa }{\partial \bxa}(\bxa)^\top\\
0&0
\end{bmatrix}
\frac{\PR(\bx) }{\qa (\bxa)^2}.
\end{multline*}
However, since we have
\begin{multline*}
\der_{ g}\PRb (\bx) =\\\begin{bmatrix}
\frac{\partial \qa }{\partial \bxa}(\bxa)^\top \frac{\Sa (\bx)}{\qa (\bxa)}
+\frac{\Sa (\bx)^\top}{\qa (\bxa)} \frac{\partial \qa }{\partial \bxa}(\bxa)
& \frac{\partial \qa }{\partial \bxa}(\bxa)^\top
	\\
\frac{\partial \qa }{\partial \bxa}(\bxa) & 0
\end{bmatrix}
\end{multline*}
and
\begin{multline*}
 \frac{\PRb (\bx) }{\qa (\bxa)^2}
\begin{bmatrix}
0&0\\
\frac{\partial \qa }{\partial \bxa}(\bxa)&0
\end{bmatrix}
+
\begin{bmatrix}
0&\frac{\partial \qa }{\partial \bxa}(\bxa)^\top\\
0&0
\end{bmatrix}
\frac{\PRb (\bx) }{\qa (\bxa)^2} =\\  \begin{bmatrix}
\frac{\partial \qa }{\partial \bxa}(\bxa)^\top \frac{\Sa (\bx)}{\qa (\bxa)}
+\frac{\Sa (\bx)^\top}{\qa (\bxa)} \frac{\partial \qa }{\partial \bxa}(\bxa)
& \frac{\partial \qa }{\partial \bxa}(\bxa)^\top
	\\
\frac{\partial \qa }{\partial \bxa}(\bxa) & 0
\end{bmatrix}
\end{multline*}
then the claim holds.
\end{proof}

\subsection{Illustrative example}

As an illustrative example, consider the case in which the vector fields $f$ and $g$ are given by
$$
f(\bx) = \begin{bmatrix}
-\bx_{a1}+\sin(\bx_{a2})\cos(\bx_{a1}) + \bx_{a2}\\
[2+\sin(\bx_{a1})]\bx_b\\
0
\end{bmatrix},g(\bx) = \begin{bmatrix}
0\\ 0\\ 1\end{bmatrix}\ .
$$
This system may be rewritten with $\bxa = (\bx_{a1},\bx_{a2})$ as 
$$
\dot {\bx}_a = \fa(\bxa) + \ga(\bxa)\bxb\ ,\ \dot {\bx}_b = u
$$
with
$$
\fa(\bxa) = \begin{bmatrix}
-\bx_{a1} + \sin(\bx_{a2})\cos(\bx_{a1}) + \bx_{a2}\\
0
\end{bmatrix}\ ,
$$
$$\ga(\bxa) = \begin{bmatrix}
0\\ 2 + \sin(\bx_{a1})
\end{bmatrix}
$$
Consider the matrix $\PRa = \begin{pmatrix}
2& 1\\1& 2
\end{pmatrix}$.
Note that if we consider
$
\Ua(\bxa) = \bx_{a1} + 2 \bx_{a2}
$,
then equation (\ref{LP241}) is satisfied with $\alphaa = \frac{1}{2 + \sin(\bx_{a1})}$.
Moreover, 
note that we have
$
v^\top \frac{\partial \Ua}{\partial \bxa}(\bxa) = 0\Leftrightarrow v_1 + 2 v_2 = 0
$.
Moreover, we have\\[0.5em]
$\displaystyle\begin{bmatrix}-2 & 1 \end{bmatrix}
\PRa \frac{\partial \fa}{\partial \bxa}(\bxa)
\begin{bmatrix}-2 \\ 1 \end{bmatrix}=-3\left[-2\frac{\partial f_{a1}}{\partial \bx_{a1}} + \frac{\partial f_{a1}}{\partial \bx_{a2}}\right]$\\[0.5em]
\null$\displaystyle\hphantom{\begin{bmatrix}-2 & 1 \end{bmatrix}
\PRa \frac{\partial \fa}{\partial \bxa}(\bxa)
\begin{bmatrix}-2 \\ 1 \end{bmatrix}}=-3.
$\\[0.5em]$\null\hfill[-2(-1 + \sin(\bx_{a2})\sin(\bx_{a1}))
 - \cos(\bx_{a1})\cos(\bx_{a2}) + 1]$\\[0.5em]
\null$\displaystyle\hphantom{\begin{bmatrix}-2 & 1 \end{bmatrix}
\PRa \frac{\partial \fa}{\partial \bxa}(\bxa)
\begin{bmatrix}-2 \\ 1 \end{bmatrix}}=-3.
$\\[0.5em]$\null\hfill[3- \sin(\bx_{a2})\sin(\bx_{a1})
 - \cos(\bx_{a1}-\bx_{a2}) ]$\\[0.5em]
\null$\displaystyle\hphantom{\begin{bmatrix}-2 & 1 \end{bmatrix}
\PRa \frac{\partial \fa}{\partial \bxa}(\bxa)
\begin{bmatrix}-2 \\ 1 \end{bmatrix}}
\leq -3$\\[0.5em]
The function  $\frac{\partial \fa}{\partial \bxa}(\bxa)$ being periodic in $\bx_{a1}$ and $\bx_{a2}$ we can assume that $\bx_{a1}$ and $\bx_{a2}$ are in a compact subset denoted $\CR_a$. This implies employing  Finsler Lemma that there exists $\rhoa$ and $\Qa$ such that inequality (\ref{LP231}) holds.
Consequently, the $\bxa$ subsystem satisfies Assumption \ref{Ass_x1Sync}.
Finally note that Assumption \ref{Ass_Pinv} is also trivially satisfied by taking $\qa(\bxa)=2+\sin(\bx_{a1})$.
From Theorem \ref{Theo_Backstepping}, it implies that there exist positive real numbers $\rhob$ and $\eta$ such that with 
$
U(\bx) = \eta (\bx_{a1} + 2 \bx_{a2}) + \frac{\bx_{b}}{2+\sin(\bx_{a1})}
$
with $\alpha(\bx) = 2+\sin(\bx_{a1})$,
equations (\ref{LP24}) and (\ref{LP23}) are satisfied.
Hence from Theorem \ref{Theo_SufCond}, the control law given in (\ref{eq_ControlLaw}) solves the local exponential synchronization problem for the $N$ identical systems that exchange information via any undirected communication graph $\mathcal G$, which is connected.

\section{Conclusion}
\label{Sec_Conclusion}
In this paper, based on recent results in \cite{AndrieuJayawardhanaPraly_TAC_TransExpStab}, we have presented necessary and sufficient conditions for the solvability of local exponential synchronization of $N$ identical affine nonlinear systems through a distributed control law. In particular, we have shown that the necessary condition is linked to the infinitesimal stabilizability of the individual system and is independent of the network topology. The existence of a symmetric covariant tensor of order two, as a result of the infinitesimal stabilizability, has allowed us to design a distributed synchronizing control law. 
When the tensor and when the controlled vector field $g$ are both constant it is shown that global exponential synchronization may be achieved.
Finally, a recursive computation of the tensor has been also discussed. 
\appendix
\if\arxiv1
\subsection{Proof of Lemma \ref{lem_Laplacian}}
\label{sec_ProofLemmaLap}
From the property of Laplacian matrix, the eigenvalues of $L$ are real and satisfy $0=\lambda_1\leq \lambda_2\leq \ldots \leq \lambda_N$. Let us take the non-zero eigenvalue $\nu>0$ of $L$ and its corresponding eigenvector $v$ in $\RR^{N}$. 
Note that we can decompose
$$
v = \begin{bmatrix}
\va \\\vb 
\end{bmatrix}\ , \ L = \begin{bmatrix}
\La & \Lb\\
\Lb^\top & \Lc
\end{bmatrix}
$$
with $\va $ and $\La$ in $\RR$. It follows that 
\begin{align}
\La \va  + \Lb\vb  &= \nu \va \\
\Lb^\top \va  + \Lc\vb  &= \nu \vb .
\end{align}
Moreover, since $\mathds{1}_N$ is an eigenvector associated to the eigenvalue $0$, 
\begin{align}
\La +\Lb\mathds{1}_{N-1} &= 0\\
\Lb^\top  + \Lc\mathds{1}_{N-1} &=0
\end{align}

Consider now a vector in $\RR^{N-1}$ defined by
$$
\tilde v = \vb  - \mathds{1}_{N-1}\va 
$$
Note that $\tilde v$ is non zero since $v$ is not colinear to $\mathds{1}_N$. By a routine algebraic computation, it follows that this vector satisfies
\\[0.5em]$
 [\Lc - \mathds{1}_{N-1} \Lb] \tilde v=  \Lc\vb  - \mathds{1}_{N-1} \Lb\vb 
$\\[0.5em]\null\hfill$
+[ \mathds{1}_{N-1}\Lb\mathds{1}_{N-1}- \Lc\mathds{1}_{N-1}]\va $
\\[0.5em]\null$\hphantom{[\Lc - \mathds{1}_{N-1} \Lb] \tilde v }=
\nu \vb  - \Lb^\top \va 
$\\[0.5em]\null\hfill$
 - \mathds{1}_{N-1} \Lb\vb +[- \mathds{1}_{N-1}\La+ \Lb]\va 
$\\[0.5em]\null$\hphantom{[\Lc - \mathds{1}_{N-1} \Lb] \tilde v }=
\nu \tilde v. 
$\\[0.5em]
This shows that $\tilde v$ is an eigenvector with the same non-zero eigenvalue of $L$. It proves the first claim of the lemma.

Note that the multiplicity of the eigenvalue $\nu$ is the same for both matrices. 
Also, if the graph is connected, then the $0$ eigenvalue of the Laplacian matrix $L$ is of multiplicity $1$ and the other eigenvalues are positive and distinct. Hence the matrix $-\bar L$ is Hurwitz. $\hfill \Box$ 

\fi
\subsection{Proof of Lemma \ref{lem_Laplacian2}}

The matrix $L$ being a balanced Laplacian matrix is positive semi-definite and its eigenvalues are real and satisfy $0=\lambda_1\leq \lambda_2\leq \ldots \leq \lambda_N$. Consequently, the principal sub-matrix $L_{2:N,2:N}$ of $L$ is also symmetric positive semi-definite (by the Cauchy's interlacing theorem). Moreover, by Kirchhoff's theorem, the matrix $L_{2:N,2:N}$, which is a minor of the Laplacian, has a determinant strictly larger than $0$ since the graph is connected. Hence, $L_{2:N,2:N}$ is positive definite. Consequently, there exists $c_1$ sufficiently small such that $A(c_1)$ is negative definite.


\if\arxiv1

\subsection{Some results from \cite{AndrieuJayawardhanaPraly_TAC_TransExpStab}}
\label{sec_Trans}
Throughout this section, we give some of the results of \cite{AndrieuJayawardhanaPraly_TAC_TransExpStab}.
Hence, we consider a system in the form
\begin{equation}
\label{eq_System_app}
\dot e = F(e,z)\ ,\quad\dot z = G(e,z)
\end{equation}
where $e$ is in $\RR^{n_e}$, $z$ is in $\RR^{n_z}$
and the functions
$F:\RR^{n_e}\times\RR^{n_z}\rightarrow \RR^{n_e}$ and
$G:\RR^{n_e}\times\RR^{n_z}\rightarrow \RR^{n_z}$ are
$C^2$. We denote by $(E(e_0,x_0,t),X(
e,z,t))$ the (unique)
solution which goes through $(e,z)$ in $\RR^{n_e}\times\RR^{n_z}$
at  $t=0$. We assume it is defined for all positive times, i.e. the
system is {\it forward complete}.

In the following, to simplify our notations, we denote by
$B_e(a)$ the open ball of radius $a$ centered at the origin in $\RR^{n_e}$.

In \cite{AndrieuJayawardhanaPraly_TAC_TransExpStab}, the following three  notions are introduced. 
\begin{list}{}{%
\parskip 0pt plus 0pt minus 0pt%
\topsep 1ex plus 0pt minus 0pt%
\parsep 0.5ex plus 0pt minus 0pt%
\partopsep 0pt plus 0pt minus 0pt%
\itemsep 1ex plus 0pt minus 0pt
\settowidth{\labelwidth}{1em}%
\setlength{\labelsep}{0.5em}%
\setlength{\leftmargin}{\labelwidth}%
\addtolength{\leftmargin}{\labelsep}%
}
\item[\textsf{
TULES-NL}]
\itshape
(Transversal uniform
local exponential stability)
\\
There exist strictly positive real numbers $r$, $k$ and
$\lambda$ such
that we have, for all $(e,x,t)$ in $\RR^{n_e}\times\RR^{n_z}\times
\RR_{\geq 0}$ with $|e|\leq r$,
\begin{equation}
\label{eq_ExpStab}
|E(e,x,t)| \leq k |e| \exp(-\lambda t)
\  .
\end{equation}
\upshape
\item[\textsf{UES-TL}]
\itshape
(Uniform exponential stability for the transversally linear system)\\
The system
\begin{equation}
\label{LP4}
\dot{\dx }= \dG(\dx)  :=  G(0,\dx)
\end{equation}
is forward complete and there exist strictly positive real numbers $\dk $  and $\tilde \lambda$
such that any solution $(\dE(\de,\dx,t),\dX(\dx,t))$  of the transversally linear system
\begin{equation}
\label{eq_System_dif}
\dot \de = \frac{\partial F}{\partial e}(0,\dx)\de\ ,\quad\dot \dx = \dG(\dx)
\end{equation}
satisfies, for all $(\de,\dx,t)$ in $\RR^{n_e}\times\RR^{n_z}\times
\RR_{\geq 0}$,
\begin{equation}\label{eq_ExpStabDrift}
|\dE(\de,\dx,t)|\leq  \dk \exp(-\tilde \lambda t)|\de|
\ .
\end{equation}
\upshape
\item[\textsf{ULMTE}]
\itshape
(Uniform Lyapunov matrix transversal equation)
\\
For all positive definite matrix $Q$, there exists a continuous
function $\PR :\RR^{n_z}\rightarrow\RR^{n_e\times n_e}$
and strictly positive real numbers $\underline{p}$ and $\overline{p}$
such that for all $\dx$ in $\RR^{n_z}$,
\begin{eqnarray}
\label{eq_TensorDerivative}
&\displaystyle
\hskip -3em
\der_{\dG} \PR(\dx) +
\PR (\dx ) \frac{\partial F}{\partial e}(0,\dx)
+ \frac{\partial F}{\partial e}(0,\dx)^\prime \PR(\dx)
\leq -Q
\\[0.5em]
\label{LP7}
&\displaystyle
\underline{p} \,  I \leq  \PR (\dx)\leq \overline{p}\,  I
\  .
\end{eqnarray}
\upshape
\end{list}

From these definitions and in the same spirit as Lyapunov second method, the following relationships have been established in \cite{AndrieuJayawardhanaPraly_TAC_TransExpStab}.
\begin{proposition}[\cite{AndrieuJayawardhanaPraly_TAC_TransExpStab},\textsf{
TULES-NL} ``$ \Rightarrow $'' \textsf{UES-TL} ]
\label{Prop_DetecNec}
 If Property \textsf{TULES-NL} holds and there exist positive
real number $c$ such that, for all $z$ in
$\RR^{n_z}$,
\begin{equation}
\label{LP1}
 \left|\frac{\partial F}{\partial e}(0,z)\right|\leq c
\ ,\quad
\left|\frac{\partial G}{\partial x}(0,x)\right|\leq c
\end{equation}
and, for all $(e,x)$ in $B_e(kr)\times\RR^{n_z}$,
\begin{equation}
\label{LP2}
\left|\frac{\partial^2 F}{\partial e\partial e}(e,z)\right|\leq c\ ,\
\left|\frac{\partial^2 F}{\partial z\partial e}(e,z)\right|\leq c\  ,\
\left|\frac{\partial G}{\partial e}(e,z)\right|\leq c
\; ,
\end{equation}
then Property \textsf{UES-TL} holds.
\end{proposition}
\vspace{0.2cm}

\begin{proposition}[\textsf{UES-TL}  ``$ \Rightarrow $'' \textsf{ULMTE}]
\label{Prop_ExistTensor}
 If  
Property \textsf{UES-TL} holds, $P$ is $C^1$  and there exists a positive real number $c$ such that
\begin{equation}
\label{LP5}
\left|\frac{\partial F}{\partial e}(0,z)\right|\leq c
\qquad \forall z\in \RR^{n_z}\ ,
\end{equation}
then Property \textsf{ULMTE} holds.
\end{proposition}
\vspace{0.2cm}

\begin{proposition}[\textsf{ULMTE} ``$\Rightarrow$'' \textsf{TULES-NL}]
\label{Prop_Lyap}
If Property \textsf{ULMTE} holds and there exist positive
real numbers $\eta  $ and $c$ such that,
for all $(e,x)$ in $B_e(\eta )\times\RR^{n_x}$,
\begin{eqnarray}
\label{LP8}
&\displaystyle \left| \frac{\partial P}{\partial x} (x)\right|\leq c
\; ,
\\[0.3em]
\label{LP6}
&\hskip -1.6 em
\displaystyle \left|\frac{\partial^2 F}{\partial e\partial
e}(e,x)\right|\leq c\; ,\
\left|\frac{\partial^2 F}{\partial x\partial e}(e,x)\right|\leq c\;  ,\
\left|\frac{\partial G}{\partial e}(e,x)\right|\leq c
\, ,\null
\end{eqnarray}
then Property \textsf{TULES-NL} holds.
\end{proposition}
\vspace{0.2cm}
\fi

\bibliographystyle{plain}

\bibliography{BibVA}

\end{document}